\documentclass[12pt]{article}

\usepackage{amsmath}
\usepackage{amssymb}

\begin{document}

\newcommand{\B}{{\mathcal B}}
\newcommand{\Bt}{{\wt{\mathcal B}}}
\newcommand{\BMW}{{\mathcal B}{\mathcal M}{\mathcal W}}
\newcommand{\SB}{{\mathcal S}{\mathcal B}}
\newcommand{\Hc}{{\mathcal H}}
\newcommand{\End}{{\rm{End}\ts}}
\newcommand{\non}{\nonumber}
\newcommand{\wt}{\widetilde}
\newcommand{\wh}{\widehat}
\newcommand{\ot}{\otimes}
\newcommand{\la}{\lambda}
\newcommand{\al}{\alpha}
\newcommand{\be}{\beta}
\newcommand{\ga}{\gamma}
\newcommand{\si}{\sigma^{}}
\newcommand{\sii}{\sigma^{-1}}
\newcommand{\pR}{{}^{\prime}\hspace{-0.15em}R}
\newcommand{\pP}{{}^{\prime}\hspace{-0.1em}P}
\newcommand{\pA}{{}^{\prime}\hspace{-0.2em}A}
\newcommand{\Pp}{P^{\ts\prime}}
\newcommand{\Rp}{R^{\ts\prime}}
\newcommand{\Ap}{A^{\prime}}
\newcommand{\ptR}{{}^{\prime}\hspace{-0.1em}{\wt R}}
\newcommand{\wR}{{\wt R}}
\newcommand{\cR}{\check R}
\newcommand{\oT}{{\overline T}}
\newcommand{\oQ}{{\overline Q}}
\newcommand{\ee}{e^{}}
\newcommand{\de}{\delta^{}}
\newcommand{\mut}{\wt{\mu}}
\newcommand{\hra}{\hookrightarrow}
\newcommand{\tpr}{t^{\tss\prime}}
\newcommand{\ve}{\varepsilon}
\newcommand{\ts}{\,}
\newcommand{\qin}{q^{-1}}
\newcommand{\tss}{\hspace{1pt}}
\newcommand{\U}{ {\rm U}}
\newcommand{\Y}{ {\rm Y}}
\newcommand{\SY}{ {\rm SY}}
\newcommand{\C}{\mathbb{C}\tss}
\newcommand{\Z}{\mathbb{Z}}
\newcommand{\ZZ}{{\rm Z}}
\newcommand{\gl}{\mathfrak{gl}}
\newcommand{\oa}{\mathfrak{o}}
\newcommand{\spa}{\mathfrak{sp}}
\newcommand{\g}{\mathfrak{g}}
\newcommand{\ka}{\mathfrak{k}}
\newcommand{\p}{\mathfrak{p}}
\newcommand{\sll}{\mathfrak{sl}}
\newcommand{\agot}{\mathfrak{a}}
\newcommand{\qdet}{ {\rm qdet}\ts}
\newcommand{\sdet}{ {\rm sdet}\ts}
\newcommand{\gr}{ {\rm gr}}
\newcommand{\sgn}{ {\rm sgn}}
\newcommand{\Sym}{\mathfrak S}

\newcommand{\Proof}{\noindent{\it Proof.}\ \ }   
\renewcommand{\theequation}{\arabic{section}.\arabic{equation}} 

\newtheorem{thm}{Theorem}[section]
\newtheorem{prop}[thm]{Proposition}
\newtheorem{cor}[thm]{Corollary}
\newtheorem{defin}[thm]{Definition}
\newtheorem{example}[thm]{Example}
\newtheorem{lem}[thm]{Lemma}

\newcommand{\bth}{\begin{thm}}
\renewcommand{\eth}{\end{thm}}
\newcommand{\bpr}{\begin{prop}}
\newcommand{\epr}{\end{prop}}
\newcommand{\ble}{\begin{lem}}
\newcommand{\ele}{\end{lem}}
\newcommand{\bco}{\begin{cor}}
\newcommand{\eco}{\end{cor}}
\newcommand{\bde}{\begin{defin}}
\newcommand{\ede}{\end{defin}}
\newcommand{\bex}{\begin{example}}
\newcommand{\eex}{\end{example}}

\newcommand{\bal}{\begin{aligned}}
\newcommand{\eal}{\end{aligned}}
\newcommand{\beq}{\begin{equation}}
\newcommand{\ben}{\begin{equation*}}

\def\beql#1{\begin{equation}\label{#1}}

\newbox\squ  
\setbox\squ=\hbox{\vrule width.3pt
            \vbox{\hrule height.3pt width.4em\kern1ex\hrule height.3pt}%
            \vrule width.3pt}
\def\endproof{%
  \ifmmode\eqno\copy\squ\smallskip
  \else{\unskip\nobreak\hfil%
    \penalty50\hskip2em\hbox{}\nobreak\hfil\copy\squ
    \parfillskip=0pt \finalhyphendemerits=0\penalty-100\smallskip}
  \fi}
\title{\Large\bf A new quantum analog of the Brauer algebra}
\author{{\sc A. I. Molev}\\[15pt]
School of Mathematics and Statistics\\
University of Sydney,
NSW 2006, Australia\\
{\tt
alexm\hspace{0.09em}@\hspace{0.1em}maths.usyd.edu.au
}
}

\date{} 
\maketitle

\vspace{12 mm}

\begin{abstract}
We introduce a new algebra $\B_l(z,q)$ depending on
two nonzero complex parameters such that $\B_l(q^n,q)$
at $q=1$ coincides with the Brauer algebra $\B_l(n)$.
We establish an analog of the Brauer--Schur--Weyl duality
where the action of the new algebra commutes with
the representation of the twisted deformation
$\U'_q(\oa_n)$ of the enveloping algebra  $\U(\oa_n)$
in the tensor power of the vector representation.
\vspace{5 mm}

\end{abstract}


\newpage

\section{Introduction}\label{sec:int}
\setcounter{equation}{0}

In the classical Schur--Weyl duality the actions of the general linear
group $GL(n)$ and the symmetric group $\Sym_l$
in the tensor power $(\C^n)^{\ot\ts l}$ are centralizers of each other.
A quantum analog of this duality is provided
by the actions of the quantized enveloping algebra $\U_q(\gl_n)$
and the Iwahori--Hecke algebra $\Hc_l(q)$ on the same space; see Jimbo~\cite{j:qu},
Leduc and Ram~\cite{lr:rh}. 

If the group $GL(n)$ is replaced with the orthogonal group $O(n)$
or symplectic group $Sp(n)$ (with even $n$)
the corresponding centralizer is generated by the action of a larger
algebra $\B_l(n)$ called the Brauer algebra originally
introduced by Brauer in \cite{b:ac}. 
The group algebra $\C[\Sym_l]$ is a natural subalgebra of $\B_l(n)$.
When $\U_q(\gl_n)$ is replaced by the quantized enveloping algebra
$\U_q(\oa_n)$ or $\U_q(\spa_n)$
then the corresponding centralizer is generated by the action of the
algebra $\BMW_l(z,q)$ introduced by Birman and Wenzl~\cite{bw:bl}
and Murakami~\cite{m:rq}; here the parameter $z$ should be appropriately specialized.
The algebraic structure and representations 
of the Brauer algebra and
the Birman--Wenzl--Murakami algebra were also studied by
Wenzl~\cite{w:sb}, Ram and Wenzl~\cite{rw:mu}, Leduc and Ram~\cite{lr:rh},
Halverson--Ram~\cite{hr:ca}, Nazarov~\cite{n:yo}.

Unlike the classical case, the algebras $\U_q(\oa_n)$ and  $\U_q(\spa_n)$
are not isomorphic to subalgebras of $\U_q(\gl_n)$. Accordingly,
the Iwahori--Hecke algebra $\Hc_l(q)$ is not a natural subalgebra of $\BMW_l(z,q)$.
On the other hand, there exist subalgebras
$\U'_q(\oa_n)$ and $\U'_q(\spa_n)$ of $\U_q(\gl_n)$ which
specialize respectively to $\U(\oa_n)$ and $\U(\spa_n)$
as $q\to 1$. They were first introduced by Gavrilik and Klimyk~\cite{gk:qd}
(orthogonal case) and by Noumi~\cite{n:ms} (both cases).
Following Noumi we call them {\it twisted quantized enveloping algebras\/}.
Although $\U'_q(\oa_n)$ and $\U'_q(\spa_n)$ are not Hopf algebras,
they are coideal subalgebras of $\U_q(\gl_n)$.
These algebras were studied in \cite{n:ms} in connection with
the theory of quantum symmetric spaces.
The algebra $\U'_q(\oa_n)$ also appears as the symmetry
algebra for the $q$-oscillator representation of
the quantized enveloping algebra $\U_q(\spa_{2m})$; 
see Noumi, Umeda and Wakayama~\cite{nuw:dp}.
Central elements and representations of $\U'_q(\oa_n)$ were studied by
Gavrilik and Iorgov~\cite{gi:ce},
Havl\'\i\v cek, Klimyk and Po\v sta~\cite{hkp:ce},
Klimyk~\cite{k:ci}.
A relationship between the algebras $\U'_q(\oa_n)$ and $\U'_q(\spa_n)$
and their affine analogs was studied in \cite{mrs:cs}.

In this paper we only consider the algebra $\U'_q(\oa_n)$ although the results
can be easily extended to the symplectic case as well.
We define a new algebra $\B_l(z,q)$ and show that its action
on the tensor power of the vector representation of $\U_q(\gl_n)$
(with $z=q^n$) commutes with that of the subalgebra 
$\U'_q(\oa_n)\subseteq\U_q(\gl_n)$.
We thus have an embedding (at least for generic $q$) of
the Iwahori--Hecke algebra $\Hc_l(q)$ into $\B_l(z,q)$.
Moreover, the algebra $\B_l(q^n,q)$ coincides
with the Brauer algebra $\B_l(n)$ for $q=1$.


This work was inspired by Arun Ram's illuminating talks in Sydney and Newcastle.
I would also like to thank him for valuable discussions.

\section{Brauer algebra and its quantum analog}\label{sec:bra}
\setcounter{equation}{0}

Let $l$ be a positive integer and $\eta$ a complex number.
An $l$-diagram $d$ is a collection
of $2l$ dots arranged into two rows with $l$ dots in each row
connected by $l$ edges such that any dot belongs to only one edge.
The product of two diagrams $d_1$ and $d_2$ is determined by
placing $d_1$ above $d_2$ and identifying the vertices
of the bottom row of $d_1$ with the corresponding
vertices in the top row of $d_2$. Let $s$ be the number of
loops in the picture. The product $d_1d_2$ is given by
$\eta^{\tss s}$ times the resulting diagram without loops.
The {\it Brauer algebra\/} $\B_l(\eta)$ is defined as the
$\C$-linear span of the $l$-diagrams with the multiplication defined above.
The dimension of the algebra is $1\cdot 3\cdots (2l-1)$.
The following presentation of $\B_l(\eta)$ is well-known; see, e.g., \cite{bw:bl}.

\bpr\label{prop:bradr}
The Brauer algebra $\B_l(\eta)$ is isomorphic to the algebra with $2l-2$ generators
$\si_1,\dots,\si_{l-1},e_1,\dots,e_{l-1}$
and the defining relations
\beql{bradr}
\begin{aligned}
\sigma_i^2&=1,\qquad e_i^2=\eta\ts \ee_i,\qquad 
\si_i\ee_i=\ee_i\si_i=\ee_i,\qquad i=1,\dots,l-1,\\
\si_i \si_j &=\si_j\si_i,\qquad \ee_i \ee_j = \ee_j \ee_i,\qquad 
\si_i \ee_j = \ee_j\si_i,\qquad
|i-j|>1,\\
\si_i\si_{i+1}\si_i&=\si_{i+1}\si_i\si_{i+1},\qquad 
\ee_i\ee_{i+1}\ee_i=\ee_i,\qquad \ee_{i+1}\ee_i\ee_{i+1}=\ee_{i+1},\\
\si_i\ee_{i+1}\ee_i&=\si_{i+1}\ee_i,\qquad \ee_{i+1}\ee_i\si_{i+1}=\ee_{i+1}\si_i,\qquad
i=1,\dots,l-2.
\end{aligned}
\non
\end{equation}
\epr

The generators $\si_i$ and $\ee_i$ correspond to the following diagrams
respectively:

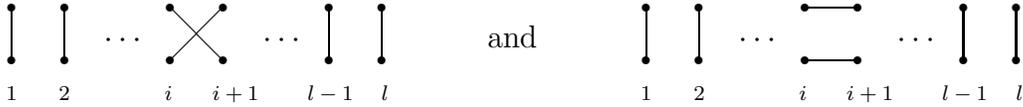
\begin{center}
\begin{picture}(400,60)
\thinlines

\put(10,20){\circle*{3}}
\put(30,20){\circle*{3}}
\put(70,20){\circle*{3}}
\put(90,20){\circle*{3}}
\put(130,20){\circle*{3}}
\put(150,20){\circle*{3}}

\put(10,40){\circle*{3}}
\put(30,40){\circle*{3}}
\put(70,40){\circle*{3}}
\put(90,40){\circle*{3}}
\put(130,40){\circle*{3}}
\put(150,40){\circle*{3}}

\put(10,20){\line(0,1){20}}
\put(30,20){\line(0,1){20}}
\put(70,20){\line(1,1){20}}
\put(90,20){\line(-1,1){20}}
\put(130,20){\line(0,1){20}}
\put(150,20){\line(0,1){20}}

\put(45,25){$\cdots$}
\put(105,25){$\cdots$}

\put(8,5){\scriptsize $1$ }
\put(28,5){\scriptsize $2$ }
\put(68,5){\scriptsize $i$ }
\put(86,5){\scriptsize $i+1$ }
\put(122,5){\scriptsize $l-1$ }
\put(150,5){\scriptsize $l$ }

\put(190,25){\text{and}}

\put(250,20){\circle*{3}}
\put(270,20){\circle*{3}}
\put(310,20){\circle*{3}}
\put(330,20){\circle*{3}}
\put(370,20){\circle*{3}}
\put(390,20){\circle*{3}}

\put(250,40){\circle*{3}}
\put(270,40){\circle*{3}}
\put(310,40){\circle*{3}}
\put(330,40){\circle*{3}}
\put(370,40){\circle*{3}}
\put(390,40){\circle*{3}}

\put(250,20){\line(0,1){20}}
\put(270,20){\line(0,1){20}}
\put(310,20){\line(1,0){20}}
\put(310,40){\line(1,0){20}}
\put(370,20){\line(0,1){20}}
\put(390,20){\line(0,1){20}}

\put(285,25){$\cdots$}
\put(345,25){$\cdots$}

\put(248,5){\scriptsize $1$ }
\put(268,5){\scriptsize $2$ }
\put(308,5){\scriptsize $i$ }
\put(326,5){\scriptsize $i+1$ }
\put(362,5){\scriptsize $l-1$ }
\put(390,5){\scriptsize $l$ }

\end{picture}
\end{center}

The subalgebra of $\B_l(\eta)$ generated by $\si_1,\dots,\si_{l-1}$
is isomorphic to the group algebra $\C[\Sym_l]$ so that $\si_i$
can be identified with the transposition $(i,i+1)$.
It is clear from the presentation that the algebra $\B_l(\eta)$
is generated by $\si_1,\dots,\si_{l-1}$ and one of the elements $\ee_i$.
The following proposition is easy to verify. We put $k=l-1$
to make the formulas more readable.

\bpr\label{prop:bradre}
The Brauer algebra $\B_l(\eta)$ is isomorphic to
the algebra with generators $\si_1,\dots,\si_{l-1},\ee_{l-1}$
and the defining relations
\beql{bradre}
\begin{aligned}
\sigma_i^2&=1,\qquad \si_i \si_j =\si_j\si_i,\qquad 
\si_i\si_{i+1}\si_i=\si_{i+1}\si_i\si_{i+1},\\
e_{k}^2&=\eta\ts \ee_{k},\qquad 
\si_{k}\ee_{k}=\ee_{k}\si_{k}=\ee_{k},\\
\ee_{k}\si_{k-1}\ee_{k}&=\ee_{k},\qquad 
\si_i \ee_{k} = \ee_{k}\si_i,\qquad i=1,\dots,k-2,\\
\ee_{k}\ts\tau\ts\ee_{k}\ts
\tau&=\tau\ts\ee_{k}
\ts\tau\ts\ee_{k},
\end{aligned}
\non
\end{equation}
where $\tau=\si_{k-1}\si_{k-2}\tss\si_{k}\si_{k-1}$. 
\epr

Note that  $\tau$
is the permutation $(k-2,k)(k-1,k+1)$ and so
the last relation of Proposition~\ref{prop:bradre} is equivalent
to the relation $\ee_{k-2}\tss\ee_{k}=\ee_{k}\tss\ee_{k-2}$
in the presentation of Proposition~\ref{prop:bradr}.

\medskip

Suppose now that $q$ and $z$ are nonzero complex numbers.

\bde\label{def:zqbra} The algebra $\B_l(z,q)$ is defined to be the algebra over $\C$
with generators $\si_1,\dots,\si_{l-1},\ee_{l-1}$
and the defining relations
\beql{zqbra}
\begin{aligned}
\sigma_i^2&=(q-\qin)\tss\si_i+1,\qquad \si_i \si_j =\si_j\si_i,\qquad 
\si_i\si_{i+1}\si_i=\si_{i+1}\si_i\si_{i+1},\\
e_{k}^2&=\frac{z-z^{-1}}{q-\qin}\ts \ee_{k},\qquad 
\si_{k}\ee_{k}=\ee_{k}\si_{k}=q\ts \ee_{k},\\
\ee_{k}\si_{k-1}\ee_{k}&=z\ts \ee_{k},\qquad 
\si_i \ee_{k} = \ee_{k}\si_i,\qquad i=1,\dots,k-2,\\
\ee_{k}\ts (zq\ts\tau^{-1}&+z^{-1}\qin\tau)\ts\ee_{k}\ts(q\ts\tau^{-1}+\qin\tau)
=(q\ts\tau^{-1}+\qin\tau)\ts\ee_{k}\ts	 (zq\ts\tau^{-1}+z^{-1}\qin\tau)\ts \ee_k,
\end{aligned}
\non
\end{equation}
where $k=l-1$ and $\tau=\si_{k-1}\si_{k-2}\tss\si_{k}\si_{k-1}$. 
\ede

We have used the same symbols for the generators of
the both algebras $\B_l(z,q)$ and $\B_l(\eta)$ since $\B_l(z,q)$ becomes
the Brauer algebra for the special case $q=1$ 
where $z$ is chosen in such a way
that the ratio $(z-z^{-1})/(q-\qin)$ takes value $\eta$ at $q=1$.
In particular, if $\eta=n$ is a positive integer then
$\B_l(q^n,q)$ coincides with $\B_l(n)$ for $q=1$.

The relations in the first line of Definition~\ref{def:zqbra}
are precisely the defining relation
of the Iwahori--Hecke algebra $\Hc_l(q)$. So we have a natural homomorphism
\beql{homhbra}
\Hc_l(q)\to \B_l(z,q).
\end{equation}
Its injectivity for generic $q$ (not a root of unity)
can be deduced from the
Schur--Weyl duality between $\U_q(\gl_n)$ and $\Hc_l(q)$; see Section~\ref{sec:qbd}.
We can therefore  regard $\Hc_l(q)$ as a subalgebra
of $\B_l(z,q)$.

The algebra $\B_l(z,q)$ has a presentation
analogous to the one given in Proposition~\ref{prop:bradr}.
However, it does not seem to exist an obvious choice of a distinguished family
of generators analogous to the $\ee_i$. As an example of such a family
we can take the elements $\ee_i$ of the algebra $\B_l(z,q)$
defined inductively by the formulas
\beql{ei}
\ee_i=\si_{i+1}\si_{i}\ts \ee_{i+1}	\sii_{i}\sii_{i+1},\qquad i=1,\dots,l-2.
\end{equation}
Then the following relations are easily deduced from Definition~\ref{def:zqbra}:
\beql{qzbraei}
\begin{aligned}
e_i^2&=\frac{z-z^{-1}}{q-\qin}\ts \ee_i,\qquad 
\si_i\ee_i=\ee_i\si_i=q\ts\ee_i,\qquad i=1,\dots,l-1,\\ 
\si_i \ee_j &= \ee_j\si_i,\qquad
|i-j|>1,\qquad\ee_i\ee_{i+1}\ee_i=\ee_i,\qquad \ee_{i+1}\ee_i\ee_{i+1}=\ee_{i+1},\\
\ee_i\si_{i+1}\ee_i&=\ee_i\si_{i-1}\ee_i=z\ts\ee_i,\qquad 
\ee_i\sii_{i+1}\ee_i=\ee_i\sii_{i-1}\ee_i=z^{-1}\ts\ee_i,
\\
\si_i\ee_{i+1}\ee_i&=zq\ts\sii_{i+1}\ee_i,\qquad 
\ee_{i+1}\ee_i\si_{i+1}=zq\ts \ee_{i+1}\sii_i.
\end{aligned}
\non
\end{equation}
The analogs of the Brauer algebra relations $\ee_i\ee_j=\ee_j\ee_i$ where $|i-j|>1$ 
have a complicated form and we shall not write them down.

\section{Quantized enveloping algebras}\label{sec:qe}
\setcounter{equation}{0}

We shall use an $R$-matrix presentation of the algebra $\U_q(\gl_n)$; see
Jimbo~\cite{j:qu} and 
Reshetikhin, Takhtajan and Faddeev~\cite{rtf:ql}.
As before, $q$ is
a nonzero complex number.
Consider the $R$-matrix
\beql{rmatrixc}
R=q\ts\sum_i E_{ii}\ot E_{ii}+\sum_{i\ne j} E_{ii}\ot E_{jj}+
(q-\qin)\sum_{i<j}E_{ij}\ot E_{ji}
\end{equation}
which is an element of $\End\C^n\ot \End\C^n$, where
the $E_{ij}$ denote the standard matrix units and the indices run over
the set $\{1,\dots,n\}$. The $R$-matrix satisfies the Yang--Baxter equation
\beql{YBEconst}
R_{12}\ts R_{13}\ts  R_{23}	=  R_{23}\ts  R_{13}\ts  R_{12},
\end{equation}
where both sides take values in $\End\C^n\ot \End\C^n\ot \End\C^n$ and
the subindices indicate the copies of $\End\C^n$, e.g.,
$R_{12}=R\ot 1$ etc.

The {\it quantized enveloping algebra\/} $\U_q(\gl_n)$ is generated
by elements $t_{ij}$ and $\bar t_{ij}$ with $1\leq i,j\leq n$
subject to the relations
\beql{defrel}
\bal
t_{ij}&=\bar t_{ji}=0, \qquad 1 \leq i<j\leq n,\\
t_{ii}\ts \bar t_{ii}&=\bar t_{ii}\ts t_{ii}=1,\qquad 1\leq i\leq n,\\
R\ts T_1T_2&=T_2T_1R,\qquad R\ts \overline T_1\overline T_2=
\overline T_2\overline T_1R,\qquad
R\ts \overline T_1T_2=T_2\overline T_1R.
\eal
\end{equation}
Here $T$ and $\overline T$ are the matrices
\beql{matrt}
T=\sum_{i,j}E_{ij}\ot t_{ij},\qquad \overline T=\sum_{i,j}
E_{ij}\ot \bar t_{ij},
\end{equation}
which are regarded as elements of the algebra $\End\C^n\ot \U_q(\gl_n)$.
Both sides of each of the $R$-matrix relations in \eqref{defrel}
are elements of $\End\C^n\ot \End\C^n\ot \U_q(\gl_n)$ and the subindices
of $T$ and $\overline T$ indicate the copies of $\End\C^n$ where
$T$ or $\overline T$ acts. 

We shall also use another $R$-matrix $\wt R$ given by
\beql{rmatrixtld}
\wt R=\qin\ts\sum_i E_{ii}\ot E_{ii}+\sum_{i\ne j} E_{ii}\ot E_{jj}+
(\qin-q)\sum_{i>j}E_{ij}\ot E_{ji}.
\end{equation}
We have the relation
\beql{relRRt}
\wt R=PR^{-1}P,
\end{equation}
where
\beql{p}
P=\sum_{i,j}E_{ij}\ot E_{ji}
\end{equation}
is the permutation operator.

The coproduct $\Delta$ on $\U_q(\gl_n)$ is defined by
the relations
\beql{copr}
\Delta(t_{ij})=\sum_{k=1}^n t_{ik}\ot t_{kj},\qquad
\Delta(\bar t_{ij})=\sum_{k=1}^n \bar t_{ik}\ot \bar t_{kj}.
\end{equation}

The universal enveloping algebra $\U(\gl_n)$ can be regarded
as a limit specialization of 
$\U_q(\gl_n)$ 
as $q\to 1$ so that
\beql{taugen1}
\frac{t_{ij}-\delta_{ij}}{q-\qin}\to E_{ij}\qquad\text{for}\quad i\geq j
\end{equation}
and
\beql{taugen2}
\frac{\bar t_{ij}-\delta_{ij}}{q-\qin}\to -E_{ij}\qquad\text{for}\quad i\leq j.
\end{equation}

Following Noumi~\cite{n:ms} introduce the {\it twisted 
quantized enveloping algebra\/} $\U'_q(\oa_n)$
as the subalgebra of $\U_q(\gl_n)$ generated by the matrix
elements $s_{ij}$ of the matrix $S=T\ts \overline T^{\ts t}$. 
More explicitly,
\beql{sij}
s_{ij}=\sum_{a=1}^n t_{ia}\ts\bar t_{ja}.
\end{equation}
It can
be easily derived from \eqref{defrel}
that the matrix $S$ satisfies the relations
\begin{align}\label{sijo}
s_{ij}&=0, \qquad 1 \leq i<j\leq n,\\
\label{sii1}
s_{ii}&=1,\qquad 1\leq i\leq n,\\
\label{rsrs}
R\ts S_1&\ts \pR\ts S_2=S_2\ts\pR S_1R,
\end{align}
where $\pR$ denotes the element obtained from $R$ by
the transposition in the first tensor factor:
\beql{rt}
\pR=q\ts\sum_i E_{ii}\ot E_{ii}+\sum_{i\ne j} E_{ii}\ot E_{jj}+
(q-\qin)\sum_{i<j}E_{ji}\ot E_{ji}.
\end{equation}
It can be shown that \eqref{sijo}--\eqref{rsrs} are defining relations
for the algebra $\U'_q(\oa_n)$; see \cite{nuw:dp}, \cite{mrs:cs}.	
A different presentation of $\U'_q(\oa_n)$
is given in the original paper by Gavrilik and Klimyk~\cite{gk:qd}.
An isomorphism between the presentations is provided by Noumi~\cite{n:ms}.

As $q\to 1$ the algebra specializes to $\U(\oa_n)$ so that
\beql{lims}
\frac{s_{ij}}{q-\qin}\to E_{ij}-E_{ji}\qquad\text{for}\quad i> j.
\end{equation}

\section{Quantum Brauer duality}\label{sec:qbd}
\setcounter{equation}{0}

We start by recalling the well-known quantum analog of the Schur--Weyl duality
between the quantized enveloping algebra	$\U_q(\gl_n)$ and
the Iwahori--Hecke algebra $\Hc_l(q)$; see \cite{j:qu}, \cite{lr:rh},
\cite[Chapter~10]{cp:gq}.
Consider the vector representation $\U_q(\gl_n)\to\End\C^n$
of the algebra	$\U_q(\gl_n)$
defined by
\beql{vect}
\begin{aligned}
t_{ii}&\mapsto \sum_{a=1}^n q^{\delta_{ia}}\ts E_{aa},\qquad
\bar t_{ii}\mapsto \sum_{a=1}^n q^{-\delta_{ia}}\ts E_{aa},\\
t_{ij}&\mapsto (q-\qin)\ts E_{ij},\qquad i>j,\\
\bar t_{ij}&\mapsto (\qin-q)\ts E_{ij},\qquad i<j.
\end{aligned}
\end{equation}
It will be convenient to interpret the representation in a matrix form.
We shall regard it as the homomorphism
\beql{vectmat}
\End\C^n\ot\U_q(\gl_n)\to\End\C^n\ot \End\C^n
\end{equation}
so that the images of the matrices $T$ and $\oT$ are given by
\beql{vectTT}
T\mapsto \pR_{01},\qquad \oT\mapsto\ptR_{01},
\end{equation}
where we label the copies of $\End\C^n$ in 
the tensor product $\End\C^n\ot \End\C^n$ by the indices $0$ and $1$,
respectively, and the left prime of the $R$-matrices denotes
the transposition in the first tensor factor; see \eqref{rt}.

Using the coproduct \eqref{copr} we consider the
representation of $\U_q(\gl_n)$ in the tensor product space
$(\C^n)^{\ot\ts l}$. In the matrix interpretation it takes the form
\beql{matmult}
\End\C^n\ot\U_q(\gl_n)\mapsto \End\C^n\ot \End\C^n\ot \cdots\ot \End\C^n,
\end{equation}
with
\beql{matmT}
T\mapsto \pR_{01}\cdots \pR_{0l},\qquad \oT\mapsto \ptR_{01}\cdots \ptR_{0l},
\end{equation}
where the multiple tensor product in \eqref{matmult} contains $l+1$ factors
labelled by $0,1,\dots,l$.

We let $\cR$ denote the element of the algebra $\End\C^n\ot\End \C^n$
defined by
\beql{rcheck}
\cR=PR=q\ts\sum_i E_{ii}\ot E_{ii}+\sum_{i\ne j} E_{ji}\ot E_{ij}+
(q-\qin)\sum_{i<j}E_{jj}\ot E_{ii}.
\end{equation}
The mapping
\beql{harep}
\si_i\mapsto \cR_{i,i+1},\qquad i=1,\dots,l-1
\end{equation}
defines a representation of the Iwahori--Hecke algebra $\Hc_l(q)$
in the space $(\C^n)^{\ot\ts l}$. 
If the parameter $q\in\C$ is generic (nonzero and not
a root of unity) then
the actions of $\U_q(\gl_n)$ and $\Hc_l(q)$ on
$(\C^n)^{\ot\ts l}$ are centralizers of each other.
Moreover, if $l<n$ then $\Hc_l(q)$ is isomorphic to
the centralizer of $\U_q(\gl_n)$.

Consider now the restriction of the $\U_q(\gl_n)$-module $(\C^n)^{\ot\ts l}$
to the subalgebra $\U'_q(\oa_n)$.	The action of the elements $s_{ij}$
is given by the formula
\beql{sact}
S\mapsto \pR_{01}\cdots \pR_{0l}\ts \wR_{0l}\cdots \wR_{01}.
\end{equation}
Clearly, the operators $\cR_{i,i+1}$ belong to the centralizer of the action of
$\U'_q(\oa_n)$. Introduce an operator $Q\in \End\C^n\ot \End\C^n$ by
\beql{q}
Q=\sum_{i,j=1}^n q^{n-2i+1}\ts E_{ij}\ot E_{ij}.
\end{equation}

\bpr\label{prop:q}
The operator $Q_{l-1,l}$ commutes with the action
of $\U'_q(\oa_n)$ on $(\C^n)^{\ot\ts l}$.
\epr

\Proof 	We need to verify that $Q_{l-1,l}$ commutes with the image of
the matrix $S$ given in \eqref{sact}.
It is sufficient to consider the case $l=2$ since  $Q_{l-1,l}$
commutes with $\pR_{0i}$ and $\wR_{0i}$ for $i\leq l-2$.
The result will follow from the relations
\beql{qrq}
Q_{12}\pR_{01}\pR_{02}\ts \wR_{02}\wR_{01}=Q_{12}=
\pR_{01}\pR_{02}\ts \wR_{02}\wR_{01}Q_{12}.
\end{equation}
Let us prove the first equality. We verify directly that
the operators $\pR$ and $\wR$ commute with each other and so
the left hand side of \eqref{qrq} equals 
$Q_{12}\pR_{01}\wR_{02}\pR_{02}\ts \wR_{01}$. We also have 
$R_{20}=\wR_{02}^{-1}$ by \eqref{relRRt} and
\beql{prr}
Q_{12}\pR_{01}=Q_{12}R_{20}.
\end{equation}
Therefore,
\beql{qr}
Q_{12}\pR_{01}\wR_{02}\pR_{02}\ts \wR_{01}=Q_{12}\pR_{02}\ts \wR_{01}=Q_{12}.
\end{equation}
The proof of the second equality in \eqref{qrq}
is similar and follows from the relations
\beql{rrq}
\pR_{02}\ts \wR_{01}Q_{12}=Q_{12}\qquad\text{and}\qquad
\pR_{01}\ts \wR_{02}Q_{12}=Q_{12}. 
\end{equation}
\endproof

The following theorem establishes a quantum analog of the Brauer
duality between the algebras $\U'_q(\oa_n)$ and $\B_l(z,q)$ at $z=q^n$.

\bth\label{thm:brad}
The mappings
\beql{sig}
\si_i\mapsto \cR_{i,i+1},\qquad i=1,\dots,l-1
\end{equation}
and
\beql{e}
\ee_{l-1}\mapsto Q_{l-1,l}
\end{equation}
define a representation of the algebra $\B_l(q^n,q)$ on the space
$(\C^n)^{\ot\ts l}$. Moreover, the actions of 
the algebras $\U'_q(\oa_n)$ and $\B_l(q^n,q)$ on this space
commute with each other.
\eth

\Proof The second statement follows from Proposition~\ref{prop:q}.
We now verify the relations of Definition~\ref{def:zqbra} with $z=q^n$
for the operators $\cR_{i,i+1}$ and $Q_{l-1,l}$.
This is done by an easy calculation for all relations
except for the last one whose proof is more involved. 
Clearly, in order
to the verify the latter we may  assume that $l=4$.
We start by proving the following
relations
\beql{qrrrrq}
Q_{34}\cR_{23}\cR_{34}\cR_{12}\cR_{23}\ts Q_{34}
=Q_{12}\ts Q_{34}+q^{n+1}\ts(q-\qin)Q_{34}(\cR_{12}+I\qin)
\end{equation}
and
\beql{qrrrrqinv}
Q_{34}\cR_{23}^{-1}\cR_{34}^{-1}\cR_{12}^{-1}\cR_{23}^{-1}\ts Q_{34}
=Q_{12}\ts Q_{34}+q^{-n-1}\ts(\qin-q)Q_{34}(\cR_{12}+I\qin),
\end{equation}
where $I$ is the identity operator.
Replacing $\cR$ by $PR$ we can rewrite
the left hand side of \eqref{qrrrrq} as
\beql{lhsl}
P_{13}P_{24}\ts Q_{12}\ts R_{14} R_{24}R_{13} R_{23}\ts Q_{34}.
\end{equation}
Introducing the diagonal matrix
\beql{d}
D=\sum_{i=1}^n q^{n-2i+1}\ts E_{ii}\in\End\C^n,
\end{equation}
we can present the operator $Q_{12}\in \End\C^n\ot\End\C^n$ as
\beql{qd}
Q_{12}=D_1 \ts\oQ_{12}=D_2 \ts\oQ_{12},
\end{equation}
where 
\beql{oq}
\oQ_{12}=\sum_{i,j=1}^n E_{ij}\ot E_{ij}=\pP_{12}=\Pp_{12},
\end{equation}
where for any operator
\beql{pap}
A=\sum_{i,j,r,s}a_{ijrs}\ts E_{ij}\ot E_{rs}\in\End\C^n\ot\End\C^n
\end{equation}
we denote
by $\pA$ and $\Ap$ the corresponding transposed operators
with respect to the first or second tensor factor:
\beql{pa}
\pA=\sum_{i,j,r,s}a_{ijrs}\ts E_{ji}\ot E_{rs}
\end{equation}
and
\beql{ap}
\Ap=\sum_{i,j,r,s}a_{ijrs}\ts E_{ij}\ot E_{sr}.
\end{equation}
We have
\beql{qr12}
Q_{12}\ts R_{14}=Q_{12}\ts \pR_{24},\qquad 
R_{23}\ts Q_{34}=D_4\ts\Rp_{24}\ts \oQ_{34}.
\end{equation}
Now the operator \eqref{lhsl} takes the form
\beql{ppqr}
P_{13}P_{24}\ts Q_{12}\ts R_{13}\ts \pR_{24} R_{24}\ts D_4\ts \Rp_{24}\ts\oQ_{34}.
\end{equation}
Multiplying directly the matrices, we get
\beql{rrdr}
\pR_{24} R_{24}\ts D_4\ts \Rp_{24}=R_{24}\ts D_4+q^{n+1}\ts (q-\qin)\ts \oQ_{24}.
\end{equation}
Furthermore,
\beql{qrd}
Q_{12}\ts R_{13}=Q_{12}\ts \pR_{23}, \qquad 
\oQ_{24}\ts Q_{34}=P_{23}\ts Q_{34},\qquad 
R_{24}\ts D_4\ts\oQ_{34}=D_3\ts \Rp_{23}\ts\oQ_{34}.
\end{equation}
The operator \eqref{ppqr} now becomes
\beql{ppqrd}
P_{13}P_{24}\ts Q_{12}\ts \pR_{23} 
\big(D_3\ts \Rp_{23}+q^{n+1}\ts (q-\qin)\ts P_{23}\big)\ts \oQ_{34}.
\end{equation}
Performing another multiplication, we write this as
\beql{ppqdq}
P_{13}P_{24}\ts Q_{12}\ts X_{23}\ts \oQ_{34}
\end{equation}
with
\beql{x}
X_{23}=D_3+q^{n+1}\ts (q-\qin)(\qin \oQ_{23}+\pR_{23}\ts P_{23}).
\end{equation}
Furthermore,
\beql{ppq}
\bal
P_{13}P_{24}\ts Q_{12}\ts X_{23}\ts \oQ_{34}&=
Q_{34}\ts P_{13}P_{24}\ts  X_{23}\ts \oQ_{34}\\
{}&=Q_{34}\ts P_{13}P_{24}\ts  X^{\ts\prime}_{24}\ts \oQ_{34}=
Q_{34}\ts P_{13}Y^{\ts\prime}_{23}\ts \oQ_{34},
\eal
\end{equation}
where $Y_{23}=P_{23}\ts  X^{\ts\prime}_{23}$.
Finally, since
\beql{px}
P_{13}Y^{\ts\prime}_{23}=Y^{\ts\prime}_{21}\ts P_{13}\qquad\text{and}\qquad
Q_{34}\ts P_{13}\ts \oQ_{34}=Q_{34},
\end{equation}
the left hand side of \eqref{qrrrrq} takes the form
$
Y^{\ts\prime}_{21}\ts Q_{34}
$
where
\beql{y}
Y^{\ts\prime}_{21}=Q_{12}+q^{n+1}\ts(q-\qin)(\cR_{12}+I\qin)
\end{equation}
thus completing the proof of \eqref{qrrrrq}.

Without giving all the details,	we note
for the proof of \eqref{qrrrrqinv} that 
$\cR_{12}^{-1}=\wR_{21}\ts P_{12}$ by \eqref{relRRt} and so the left hand side
of the relation can be written as
\beql{qrrrqpp}
Q_{34}\ts \wR_{32} \wR_{42} \wR_{31}  \wR_{41}  \ts Q_{12}\ts P_{13}P_{24}.
\end{equation}
Then we proceed in the same manner as for the proof of \eqref{qrrrrq}
modifying the argument appropriately.

Now, \eqref{qrrrrq} and \eqref{qrrrrqinv} give
\beql{qrriq}
Q_{34}\big(q^{-n-1}\cR_{23}\cR_{34}\cR_{12}\cR_{23}+
q^{n+1}\cR_{23}^{-1}\cR_{34}^{-1}\cR_{12}^{-1}\cR_{23}^{-1}\big)\ts Q_{34}
=\big(q^{-n-1}+	q^{n+1}\big) \ts Q_{12}\ts Q_{34}.
\non
\end{equation}
The proof of the theorem will be completed by verifying the following
relations
\beql{qqrrrr}
Q_{12}Q_{34}\big(q^{-1}\cR_{23}\cR_{34}\cR_{12}\cR_{23}+
q\ts\cR_{23}^{-1}\cR_{34}^{-1}\cR_{12}^{-1}\cR_{23}^{-1}\big)
=\big(q^{-3}+	q^{3}\big) \ts Q_{12}\ts Q_{34}
\end{equation}
and
\beql{rrrrqq}
\big(q^{-1}\cR_{23}\cR_{34}\cR_{12}\cR_{23}+
q\ts\cR_{23}^{-1}\cR_{34}^{-1}\cR_{12}^{-1}\cR_{23}^{-1}\big)Q_{12}Q_{34}
=\big(q^{-3}+	q^{3}\big) \ts Q_{12}\ts Q_{34}.
\end{equation}
The arguments are similar for both relations, so we only give
a proof of \eqref{qqrrrr}. Since  $Q_{12}Q_{34}=D_1 D_3\ts \oQ_{12}\oQ_{34}$
we may replace $Q_{12}Q_{34}$ with $\oQ_{12}\oQ_{34}$ in \eqref{qqrrrr}.
We have
\beql{oqoqrrrr}
\oQ_{12}\oQ_{34}\cR_{23}\cR_{34}\cR_{12}\cR_{23}
= \oQ_{12}\oQ_{34}P_{13}P_{24} R_{14} R_{24} R_{13} R_{23}.
\end{equation}
Furthermore,
\beql{qqpp}
\bal
\oQ_{12}\oQ_{34}P_{13}P_{24}&=\oQ_{12}P_{13}P_{24}\oQ_{12}=
\oQ_{12}\oQ_{23}P_{24}\oQ_{12}\\
{}&=\oQ_{12}P_{24}\oQ_{34}\oQ_{12}
= \oQ_{12}P_{24}\oQ_{12}\oQ_{34}=\oQ_{12}\oQ_{34}.
\eal
\end{equation}
The right hand side of \eqref{oqoqrrrr} now takes the form
$
\oQ_{12}\oQ_{34}R_{14} R_{24} R_{13} R_{23}.
$
Using the relations
\beql{oqr}
\oQ_{12}R_{14}=\oQ_{12}\pR_{24},\qquad \oQ_{12}\oQ_{34}R_{13}=\oQ_{12}\oQ_{34}\Rp_{14}
=\oQ_{12}\oQ_{34}\pR^{\ts\prime}_{24},
\end{equation}
we can write it as
\beql{qqprrp}
\oQ_{12}\oQ_{34}V_{24}R_{23}=\oQ_{12}\oQ_{34}V^{\ts\prime}_{23}R_{23},
\end{equation}
where
\beql{vpr}
V_{23}=\pR^{\ts\prime}_{23} \pR_{23} R_{23}.
\end{equation}

Now we transform the second summand on the left hand side of	\eqref{qqrrrr}
in a similar manner. Since $\cR^{-1}=P\wR$ we have
\beql{qqrinv}
\oQ_{12}\oQ_{34}\ts\cR_{23}^{-1}\cR_{34}^{-1}\cR_{12}^{-1}\cR_{23}^{-1}
= \oQ_{12}\oQ_{34}\wR_{14} \wR_{24} \wR_{13} \wR_{23}.
\end{equation}
Exactly as above, we write this operator in the form
\beql{qqprrw}
\oQ_{12}\oQ_{34}W^{\ts\prime}_{23}\wR_{23},
\end{equation}
where
\beql{wpr}
W_{23}=\ptR^{\ts\prime}_{23} \ptR_{23} \wR_{23}.
\end{equation}
A direct calculation shows that
\beql{wv}
\qin V^{\ts\prime}R+q\ts W^{\ts\prime}\wR=\big(q^{-3}+	q^{3}\big)\ts I,
\end{equation}
completing the proof. \endproof


\end{document}